\documentclass[preprint,12pt]{elsarticle}

\usepackage{amsthm}
\usepackage{epsfig, graphicx, xcolor}
\usepackage{amssymb}
\usepackage{amsmath, amsfonts}

\def\P{{\rm \hbox{I\kern-.2em\hbox{P}}}}
\def\H{{\rm \hbox{I\kern-.2em\hbox{H}}}}
\def\R{{\rm \hbox{I\kern-.2em\hbox{R}}}}
\def\N{{\rm \hbox{I\kern-.2em\hbox{N}}}}
\def\Z{{\rm {{\rm Z}\kern-.28em{\rm Z}}}}

\newcommand{\be}{\begin{equation}}
\newcommand{\ee}{\end{equation}}

\def\sep{ : \;\;}
\newcommand {\RR} {\mathbb R}
\newcommand {\NN} {\mathbb N}

\newcommand{\K}{{\bf k}}

\newtheorem {theorem} {Theorem}
\newtheorem {remark} {Remark}
\newtheorem {lemma} {Lemma}

\newtheorem {problem} {Problem}
\newtheorem {definition} {Definition}

\newcommand {\ttbox} [1] {\mbox {\ttfamily*1}} 

\newtheorem{assumption}{Assumption}

\journal{Journal of Approximation Theory}

\begin{document}

\begin{frontmatter}

\title{Kolmogorov's problem for completely and multiply monotone functions and the Markov moment problem}

\author[Babenko V.]{Vladyslav~Babenko}

\author[Babenko Y.]{Yuliya~Babenko\corref{cor1}  }

\author[Kovalenko]{Oleg~Kovalenko}

\cortext[cor1]{Corresponding author}

\address[Babenko V.]{Department of Mathematics and Mechanics, Dnepropetrovsk National University, Gagarina pr., 72, Dnepropetrovsk, 49010, UKRAINE}
\address[Babenko Y.]{Department of Mathematics, Kennesaw State University, 1100 South Marietta Pkwy, MD \# 9085, Marietta, GA, 30060, USA}
\address[Kovalenko]{Department of Mathematics and Mechanics, Dnepropetrovsk National University, Gagarina pr., 72, Dnepropetrovsk, 49010, UKRAINE}


	\begin{abstract} In this paper, we present the solution to Kolmogorov's problem for the classes of multiply monotone and completely monotone functions together with its connections to the Markov moment problem, Hermite-Birkhoff interpolation problem, and other extremal problems.
	\end{abstract}

\end{frontmatter}

\section{Introduction}

Many problems in Analysis and Approximation Theory are related to studying the necessary and sufficient conditions that guarantee existence of a function from a given class with a prescribed set of characteristics, and to studying the properties of the set of all functions with these characteristics. 
	Examples of such problems include, but are not limited to, interpolation of prescribed values at given points by functions from certain class (for example by polynomials, splines, perfect splines, etc.), Hermite and Birkhoff interpolation, various moment problems, the Nevanlinna-Pick problem, and many others. 
		
	Kolmogorov's problem about necessary and sufficient conditions to guarantee the existence of a function from a given class, which has prescribed values of norms of derivatives of given orders, can be considered as another problem of this type. It is possible that this question was motivated by Kolmogorov's work on the problem of dependence between norms of consecutive derivatives (problem on inequalities between norms of derivatives), that first appeared in the works of Hadamard, Hardy and Littlewood, and Landau. 
	
In this paper, we present the solution to Kolmogorov's problem for the classes of multiply monotone and completely monotone functions and provide its connections and applications to other extremal problems of Analysis, including the Markov moment problem, Hermite-Birkhoff interpolation problem, and others.

The paper is organized as follows. In Section \ref{S2} we introduce some necessary notation and present the statement of Kolmogorov's problem together with its alternative forms. Section \ref{S3} contains auxiliary results related to the Markov moment problem. Sections \ref{S4} and \ref{S5} introduce the classes of completely monotone and multiply monotone function, respectively, and present integral representations of the functions from these classes. In Section \ref{S6} we present the connection between the two classes of functions, which allows us to transfer the solution of Kolmogorov's problem from one class to the other. Section \ref{S7} studies the extremal properties of splines that provide the existence of the functions with prescribed properties. In Section \ref{S8} we provide the solution to Kolmogorov's problem for classes of completely monotone and multiply monotone functions. Section \ref{S9} presents some application of the obtained results, namely, to Hermite-Birkhoff interpolation problem, the problem of sharp estimates for intermediate moments (which belongs to the class of problems of finding extremal values of integrals), and extremal distribution functions.

\section{Statement of Kolmogorov's problem} \label{S2}
For $r\in\NN$ we let $L^r_{\infty,\infty}(\RR_-):=L_\infty(\RR_-)\bigcap L^r_{\infty}(\RR_-),$ $\|\cdot\|:=\|\cdot\|_{L_\infty(\RR_-)}$.
We consider Kolmogorov's problem stated as follows:

\begin{problem}[Kolmogorov's problem]
	Let the class of functions $X\subset L^r_{\infty,\infty}(\RR_-)$ and a set of $d$ integers
	$0\leq k_1 < k_2 < \ldots < k_{d}\leq r$ be given.
	The problem is to find necessary and sufficient conditions on the set of positive real numbers $M_{k_1}, M_{k_2}, \dots, M_{k_{d}}$
	in order to guarantee the existence of a function $x \in X$ such that $$\|x^{(k_i)}\|=M_{k_i},\qquad i=1,\dots, d.$$
\end{problem}

For $d\in\NN$ and integers $0\leq k_1<k_2<\ldots<k_{d}\leq r$ we set ${\bf k}:=(k_1,\ldots, k_{d})$ and $M_{{\bf k}} := \{M_{k_1},\ldots ,M_{k_{d}}\}.$

In addition, for a given function $x\in X$ we define $$M_{\bf k}(x):=\left (M_{k_1}(x),\ldots ,M_{k_d}(x)\right ),$$ where $$M_{k_i}(x)=\| x^{(k_i)}\|,\qquad i=1,\ldots, d.$$

\begin{definition}
The set $M_{{\bf k}}$ is called {admissible} for the class $X\subset L^r_{\infty,\infty}(\RR_-)$, if there exists a function $x\in X$  such that $\left\|x^{(k_i)}\right\|=M_{k_i}$, $i=1,2,\dots,d$
(or, for short, $M_{{\bf k}}(x)=M_{{\bf k}}$). 
\end{definition}

By $A_{\K}(X)$ we denote the collection of all admissible for the class $X$ sets $M_{{\bf k}}$.

Using the above notation, Kolmogorov's problem can be reformulated as follows. 
\begin{problem}
For the given class of functions $X\subset L^r_{\infty,\infty}(\RR_-)$ and an arbitrary set of $d$ integers $\K$, the problem is to characterize the set $A_{\K}(X)$.
\end{problem}

In addition, we consider Kolmogorov's problem in alternative form.
\begin{problem}[Kolmogorov's problem, alternative form]
For the given class of functions $$X\subset L^r_{\infty,\infty}(\RR_-)$$ and an arbitrary system of $d$ integers $\K$, the problem is to find ``minimal'' set $F_{\K}(X)\subset X$ such that $$A_{\K}(X) = \left\{M_{\K}(x)\colon x\in F_{\K}(X)\right\}.$$
\end{problem}

The history of the problem as well as the list of cases, when the results are known, can be found, for instance, in~\cite{Kovalenko13a}.

\section{Auxiliary results}\label{S3}

In this section, we present some useful notation and results related to the Markov problem. The results are classical and their proofs can be found, for example, in \cite{Karlin66,Krein73,Akhiezer65}. 


\begin{definition}\label{tchebishovSystem}
	The system of functions $u_1,\dots, u_n$ is called {Chebyshev system} on $\RR_+:=[0,\infty)$, if the functions $u_1,\dots, u_n$ are continuous on $\RR_+$ and 	$$\det \left[\|u_i(t_j)\|_{i,j=1}^n\right]>0$$ for an arbitrary set of points $0\leq t_1<\ldots<t_n<\infty$.
\end{definition}

\begin{problem}[the Markov moment problem]\label{markovProblem} 
	 Let a Chebyshev system of functions $u_1,\dots,u_n$ on $\RR_+$ be given. The problem is to find necessary and sufficient conditions on the set of numbers ${\bf c} = (c_1,\dots,c_n)\in\RR^{n}$, in order to guarantee the existence of a function $\sigma$ from the set $\mathcal{D}$ of nonnegative, nondecreasing functions of bounded variation  so that
	\begin{equation}\label{moments}
		c_k = \int\limits_0^\infty u_k(t)d\sigma(t),\qquad k=1,\dots,n.
	\end{equation}
\end{problem}

\begin{remark} 
	Note that in the statement of Problem~\ref{markovProblem} only those functions $\sigma\in\mathcal{D}$, for which integrals in~\eqref{moments} are absolutely convergent, are considered.\end{remark}

\begin{definition}\label{momentSpace} 
	By $$\mathcal{M}_n = \mathcal{M}(u_1,\dots,u_n)$$ we denote the set of all points $(c_1,\dots,c_n)\in\RR^{n}$,
	such that there exists a function $\sigma\in\mathcal{D}$ for which conditions~\eqref{moments} are satisfied. The set $\mathcal{M}_n$ is called the {\bf moment space}. \end{definition}

For a given system of continuous on $\RR_+$ functions $\left\{u_k\right\}_{k=1}^n$, we set 
$$\mathcal{C}_n = \mathcal{C}(u_1,\dots,u_n):=\left\{(u_1(t),\dots,u_n(t)),\, t\in \RR_+\right\}.$$

Note that the set $\mathcal{C}_n$ is a subset of moment space $\mathcal{M}(u_1,\dots,u_n)$, which is generated by Heaviside step functions $$H_t(x):=\begin{cases} 1, & x\geq t \\ 0, & x < t\end{cases},$$ where $x,t\geq 0$.

The following theorem can be found in~\cite[Chapter 5, \S~2]{Karlin66}.

\begin{theorem}\label{conicalcone}
	Assume that there is a system of continuous function $\left\{u_k\right\}_{k=1}^n$ defined on $\RR_+$  (not necessarily Chebyshev system). 
	Then the moment space $\mathcal{M}(u_1,\dots,u_n)$ coincides with the convex cone generated by
	$\mathcal{C}_n$.
\end{theorem}

From Theorem ~\ref{conicalcone} and Caratheodory's theorem it follows that every element of the set $\mathcal{M}(u_1,\dots,u_n)$ can be represented as a linear combination with positive coefficients of no more than $n+1$ points of the curve $\mathcal{C}(u_1,\dots,u_n)$.

Therefore, an arbitrary point ${\bf c}\in \mathcal{M}_n$ can be represented in the following form:
\begin{equation}\label{representation}
	{\bf c} = \sum\limits_{k=1}^ma_ku(t_k),\qquad m\leq n+1,
\end{equation}
where $0\leq t_1<\ldots<t_m$, $a_k>0$, $k=1,\dots,m$, and $u(t) := (u_1(t), u_2(t),\dots,u_n(t))$.

We also need the following definitions.
\begin{definition}
	Numbers $t_k$, $k=1,\dots,m$, from representation ~\eqref{representation} are called the roots of representation~\eqref{representation}.
\end{definition}

\begin{definition}\label{index} 
	Index $I({\bf c})$ of the point ${\bf c}\in \mathcal{M}_n$ is the minimal number of elements from $\mathcal{C}_n$, that need to be used to represent 	${\bf c}$ as a convex combination.  Note that the points $(u_1(t),\dots,u_n(t))$ with $t > 0$ are counted as $1$, and with $t = 0$ are counted as $\frac12$.
\end{definition}

\begin{assumption}\label{uConditions}
	We assume that the functions $\left\{u_k\right\}_{k=1}^n$ satisfy the following three conditions.
	\begin{enumerate}
		\item\label{uCondition1} Both systems $\left\{u_k\right\}_{k=1}^n$ and $\left\{u_k\right\}_{k=1}^{n-1}$ are Chebyshev systems on $\RR_+$.
		\item There exists a polynomial $u(t)=\sum\limits_{k=1}^na_ku_k(t)$, which satisfies conditions $u(t)>0$, $t\in \RR_+$ , and 
			$\varliminf\limits_{t\to\infty} u(t) > 0$.
		\item $\lim\limits_{t\to\infty}\frac{u_k(t)}{u_n(t)} = 0$, $k=1,\dots,n-1$.
	\end{enumerate}
\end{assumption}

If the functions $u_1,\dots, u_n$ satisfy all three conditions of Assumption~\ref{uConditions}, then the following theorems hold (these results are contained in~\cite[Chapter 5, \S~4]{Karlin66}; for simplicity of references we formulate them as theorems).
\begin{theorem}\label{momentSpaceBoundary}
	Nonzero vector ${\bf c}$ belongs to the boundary $\partial\mathcal{M}_n \bigcap \mathcal{M}_n$ of the moment space $\mathcal{M}_n$ if and only if $I({\bf c}) < \frac n 2$.
\end{theorem}

\begin{theorem}\label{momentSpaceMainRepresentation}
	 Let ${\bf c}\in {\rm int}\,\mathcal{M}_n$. Then there exists a representation of ${\bf c}$ with index $\frac n 2$.
\end{theorem}

\begin{theorem}\label{momentSpaceCanonicalRepresentation}
	Let ${\bf c}\in {\rm int}\,\mathcal{M}_n$. Then for all $t^*>0$ there exists a representation of vector ${\bf c}$ with index $\frac {n + 1} 2$,
	that has $t^*$ as a root.
\end{theorem}

We need the following lemma.

\begin{lemma}\label{powerFunctions}
Let $d\in\NN$ and integer numbers $0= k_1<k_2<\ldots<k_d$ be given. Then the system of functions $\left\{u_i\right\}_{i=1}^d$, $u_i(t) = t^{k_i}$, $i=1,\dots, d$, satisfies the conditions of assumption~\ref{uConditions}. 
\end{lemma}

It is easy to see that second and third conditions of Assumption~\ref{uConditions} hold. In order to prove the first condition we need to prove that the following inequality 
\begin{equation}\label{ourDetInequality}
\begin{vmatrix}
    1 & t_1^{k_2} & t_1^{k_3} & \dots  & t_1^{k_d} \\
    1 & t_2^{k_2} & t_2^{k_3} & \dots  & t_2^{k_d} \\
    \vdots & \vdots & \vdots & \ddots & \vdots \\
    1 & t_d^{k_2} & t_d^{k_3} & \dots  & t_d^{k_d} \\
\end{vmatrix}
> 0
\end{equation}
holds for all $0\leq t_1<t_2<\ldots<t_d$. 
In \cite[Chapter~1, \S 3]{Karlin66} it is proved that for all $0<t_1<\ldots<t_d$ and all $-\infty<\alpha_1<\ldots < \alpha_d< \infty$ 
\begin{equation}\label{detInequality}
\begin{vmatrix}
    t_1^{\alpha_1} & t_1^{\alpha_2} & t_1^{\alpha_3} & \dots  & t_1^{\alpha_d} \\
    t_2^{\alpha_1} & t_2^{\alpha_2} & t_2^{\alpha_3} & \dots  & t_2^{\alpha_d} \\
    \vdots & \vdots & \vdots & \ddots & \vdots \\
    t_d^{\alpha_1} & t_d^{\alpha_2} & t_d^{\alpha_3} & \dots  & t_d^{\alpha_d} \\
\end{vmatrix}
> 0
\end{equation}

In the case when $t_1>0$, inequality~\eqref{ourDetInequality} follows from~\eqref{detInequality}. If $t_1=0$ then 
$$
\begin{vmatrix}
    1 & t_1^{k_2} & t_1^{k_3} & \dots  & t_1^{k_d} \\
    1 & t_2^{k_2} & t_2^{k_3} & \dots  & t_2^{k_d} \\
    \vdots & \vdots & \vdots & \ddots & \vdots \\
    1 & t_d^{k_2} & t_d^{k_3} & \dots  & t_d^{k_d} \\
\end{vmatrix}
= 
\begin{vmatrix}
    1 & 0 & 0 & \dots  & 0 \\
    1 & t_2^{k_2} & t_2^{k_3} & \dots  & t_2^{k_d} \\
    \vdots & \vdots & \vdots & \ddots & \vdots \\
    1 & t_d^{k_2} & t_d^{k_3} & \dots  & t_d^{k_d} \\
\end{vmatrix}
=
\begin{vmatrix}
    t_2^{k_2} & t_2^{k_3} & \dots  & t_2^{k_d} \\
    \vdots & \vdots & \ddots & \vdots \\
    t_d^{k_2} & t_d^{k_3} & \dots  & t_d^{k_d} \\
\end{vmatrix}
>0
$$
(the last inequality follows from~\eqref{detInequality}). The lemma is proved. $\square$
\begin{remark}
	Let $d\in\NN$ and integers $0\leq k_1<k_2<\ldots<k_d$ be given. From now on we consider only systems of functions $\left\{u_i\right\}_{i=1}^d$ of the following form: $$u_i(t) = t^{k_i}, \qquad i=1,\dots, d.$$
	In the case when $k_1=0$, such systems, by Lemma~\ref{powerFunctions}, satisfy the conditions of Assumption~\ref{uConditions}, and, hence,Theorems~\ref{momentSpaceBoundary}, \ref{momentSpaceMainRepresentation}, and~\ref{momentSpaceCanonicalRepresentation} hold for them.
\end{remark}
	
\section{Completely monotone functions and their integral representation}\label{S4}
\begin{definition} 
	Infinitely many times differentiable on $\RR_-$ function is called completely (or, sometimes, absolutely) monotone, if all its derivatives are nonnegative on $\RR_-$. 
	By $CM(\RR_-)$ we denote the class of completely monotone on $\RR_-$ functions.
\end{definition}

The following useful integral representation of completely monotone functions was proved by Bernstein~\cite{Bernshtein}. 

\begin{theorem}\label{AMrepresentation}
	A funciton $x(t)$ is completely monotone on  $\RR_-$ if and only if it admits the following representation
	\begin{equation}\label{Bern}
		x(t)= \int_0^{\infty} e^{t u} d\beta(u),\;\;\; t\in\RR_-,
	\end{equation}
	where $\beta(u)$ is nondecreasing bounded function.
\end{theorem}

Note that, by definition, uniform norms of a completely monotone on $\RR_-$ function $x(t)$ and all its derivatives are achieved at zero. Due to the integral representation~\eqref{Bern} it implies that for $k=0,1,\dots$ 
$$\|x^{(k)}\| = x^{(k)}(0) =  \int_0^{\infty} u^k d\beta(u).$$

Therefore, we have the following theorem, which shows the connection between a solution to Kolmogorov's problem on the class of completely monotone functions and the Markov moment problem.
\begin{theorem}\label{AMandMomentProblem}
	Let $d\in\NN$, $\K = (0\leq k_1 < k_2<\ldots<k_d)$, be given. Then the set $A_{\K}(CM(\RR_-))$ of admissible sets for the class $CM(\RR_-)$ coincides with the moment space $\mathcal{M}(t^{k_1},\dots,t^{k_d})$.
\end{theorem}

Note that the Heaviside function $H_a$, $a > 0$, generates completely monotone function $e^{at}$. For convenience, instead, we consider completely monotone function $e^{\frac t a}$, which is generated by the Heaviside function $H_{a^{-1}}$, $a>0$. 

\begin{definition}
	A function of the form $$\phi(CM(\RR_-),{\bf a}, {\bf \lambda};t):=\sum\limits_{s=1}^m \lambda_s a_s^r e^{a_s^{-1} t},$$ where
	$\lambda_s, a_s > 0$, $s=1,\dots,m$, ${\bf \lambda} = (\lambda_1,\dots,\lambda_m)$, ${\bf a} = (a_1,\dots,a_m)$ 
	is called $CM(\RR_-)$-perfect spline of order $r\in\NN$ with $m$ knots $-a_1,\dots,-a_m$.
	
	If $\phi(t)$ is $CM(\RR_-)$-perfect spline with $m$ knots and $C > 0$, then function $C+\phi(t)$ is called
	$CM(\RR_-)$-perfect spline with $m+\frac 1 2$ knots.
\end{definition}

\section{Multiply monotone functions and their integral representation}\label{S5}

By $L_{\infty,\infty}^{r,\smallsmile}(\RR_-)$ we denote the class of functions $x\in L_{\infty,\infty}^{r}(\RR_-)$ such that for $k=0,\dots, r-1$ the derivatives $x^{(k)}$ are nondecreasing and concave (see~\cite{Williamson}). The functions from this class are called {\it {multiply monotone}} and were defined in such form by Williamson ~\cite{Williamson}.


He proved the following theorem (see, \cite{Williamson}), which provides an integral representation of functions from this class. 
\begin{theorem}
	A function $y(t)$ belongs to $L_{\infty,\infty}^{r,\smallsmile}(\RR_-)$ if and only if
	\begin{equation}\label{Will}
		y(t)=\frac{1}{r!}\displaystyle \int_0^{\infty} \left[(1+ut)_+\right]^rd\beta(u),\qquad t\in\RR_-,
	\end{equation}
	where $\beta(u)$ is nondecreasing bounded function.
\end{theorem}

Note that the Heaviside function $H_a$, $a > 0$, generates  $\frac 1{r!}(1+at)^r_+$. It will be convenient for us to use the function $\frac 1{a^rr!}(a+t)^r_+$, which is generated by $H_{a^{-1}}$, $a>0$.

\begin{definition}
	Function $$\phi(L_{\infty,\infty}^{r,\smallsmile}(\RR_-),{\bf a}, {\bf \lambda};t):=\frac 1 {r!}\sum\limits_{s=1}^m \lambda_s (a_s + t)^r_+,$$ where 
	$\lambda_s, a_s > 0$, $s=1,\dots,m$, ${\bf \lambda} = (\lambda_1,\dots,\lambda_m)$, ${\bf a} = (a_1,\dots,a_m)$, is called $L_{\infty,\infty}^{r,\smallsmile}(\RR_-)$-perfect spline of order $r\in\NN$ with $m$ knots $-a_1,\dots,-a_m$.
	
	If $\phi(t)$ is a $L_{\infty,\infty}^{r,\smallsmile}(\RR_-)$-perfect spline with $m$ knots and $C > 0$, then the function $C+\phi(t)$ is called a $L_{\infty,\infty}^{r,\smallsmile}(\RR_-)$-perfect spline with $m+\frac 1 2$ knots.
\end{definition}

\section{Connection between classes $L_{\infty,\infty}^{r,\smallsmile}(\RR_-)$ and $CM(\RR_-)$} \label{S6}

For numbers $a_1,\dots, a_d$ by ${\rm diag}(a_1,\dots,a_d)$ we denote the square diagonal matrix of order $d$ with $a_1,\dots, a_d$ on its main diagonal. For the given vector ${\bf c} \in \RR^d$ by ${\rm diag}(a_1,\dots, a_d)\,{\bf c}$ we denote the result of multiplication of the matrix ${\rm diag}(a_1,\dots,a_d)$ by a vector column ${\bf c}$. For a set $A\subset \RR^d$ we let $${\rm diag}(a_1,\dots, a_d)\,A:=\left\{{\rm diag}(a_1,\dots, a_d)\,{\bf c}\colon {\bf c} \in A\right\}.$$

Taking into account representations \eqref{Bern}, \eqref{Will}, and the fact  that norms of all the derivatives of orders $k_1,\dots,k_d$ for the functions from the classes $L_{\infty,\infty}^{r,\smallsmile}(\RR_-)$ and $CM(\RR_-)$ are achieved at zero, we obtain the following theorem.

\begin{theorem}\label{connection}
	Let integers $0\leq k_1<k_2<\ldots<k_d\leq r$, $\K = (k_1,\dots,k_d)$, be given. Then 
	$$A_{\K}(CM(\RR_-)) = {\rm diag}((r-k_1)!,\dots, (r-k_d)!) A_{\K}(L_{\infty,\infty}^{r,\smallsmile}(\RR_-)).$$
\end{theorem}

\begin{remark}\label{corresp}
	Let vectors ${\bf a}, {\bf \lambda}$ of positive numbers be given. Then 
	$$M_{\K}(\phi(CM(\RR_-),{\bf a},{\bf \lambda};t)) = {\rm diag}((r-k_1)!,\dots, (r-k_d)!) M_{\K}(\phi(L_{\infty,\infty}^{r,\smallsmile}(\RR_-),{\bf a},{\bf \lambda};t)).$$
\end{remark}

%

From now on, by $X$ we denote either $L_{\infty,\infty}^{r,\smallsmile}(\RR_-)$ or $CM(\RR_-)$.

\section{Some properties of $X$-perfect splines} \label{S7}

The next lemma (comparison type result) establishes some extremal properties of perfect splines introduced in Sections \ref{S4} and \ref{S5}.

\begin{lemma}\label{l::1}
	Let $x\in {X}$, $0 \leq s < k_1<\ldots<k_{2n}\leq r$. Let ${X}$-spline $\phi(t)$ with no more than $n$ knots be such that 
\begin{equation}\label{equalNorms}
\left\|\phi^{(k_i)}\right\| = \left\|x^{(k_i)}\right\|,\qquad i = 1,\dots,2n.
\end{equation}
Then
	\begin{equation}\label{l1}
		\left\|\phi^{(s)}\right\|\leq \left\|x^{(s)}\right\|
	\end{equation}
	and if $k_{2n} <r$, then
	\begin{equation}\label{l2}
		\left\|\phi^{(r)}\right\|\leq \left\|x^{(r)}\right\|.
	\end{equation}	
	If the equality in inequality~\eqref{l1} holds, then $x^{(s)} \equiv \phi^{(s)}$. If $k_{2n}<r$ and equality in inequality~\eqref{l2} holds, or the number of knots of the spline $\phi$ is less than $n$, then $x^{(k_1)} \equiv \phi^{(k_1)}$.
\end{lemma}
{\bf Proof.} Due to Theorem~\ref{connection}, it is sufficient to prove the statement of the theorem in the case when $X=L_{\infty,\infty}^{r,\smallsmile}(\RR_-)$. First we assume that $k_{2n}=r$. Assume to the contrary that $x^{(s)} \ne \phi^{(s)}$ and $\|x^{(s)}\| \leq \|\phi^{(s)}\|$.
Set $\Delta(t):=x(t)-\phi(t)$. In order to obtain a contradiction we count the number of sign changes of the difference $\Delta(t)$ and its derivatives.

First of all, we observe that by the definition of the spline $\phi$, we have $\phi^{(s)}(-a_1)=0$ (where $-a_1$ is the leftmost knot of the spline $\phi$). Besides that, $x^{(s)}(-a_1)\geq 0$, and hence $\Delta^{(s)}(-a_1)\geq 0$. 

By assumption
$$\Delta^{(s)}(0)=x^{(s)}(0)-\phi^{(s)}(0)=\left\|x^{(s)}\right\|-\left\|\phi^{(s)}\right\|\leq 0.$$ 
These observations imply that there exists a point $t_{s+1}^1\in(-a_1,0)$ such that $\Delta^{(s+1)}(t_{s+1}^1)< 0$. Besides that, 
$\Delta^{(s+1)}(-a_1)\geq 0$. Therefore, there exists a point $t_{s+2}^1\in(-a_1,0)$ such that $\Delta^{(s+2)}(t_{s+2}^1)<0$. Repeating the same argument, we obtain that there exists a point $t_{k_1}^1\in(-a_1,0)$ such that $\Delta^{(k_1)}(t_{k_1}^1)< 0$. In addition, $\Delta^{(k_1)}(-a_1)\geq
0$ and $\Delta^{(k_1)}(0)= 0$ by Lemma's assumptions. Hence, there are two points $-a_1<t^1_{k_1+1}<t^2_{k_1+1}<0$ such that $\Delta^{(k_1+1)}(t^1_{k_1+1})<0$ and $\Delta^{(k_1+1)}(t^2_{k_1+1})>0$. This sign distribution will remain up to the level $k_2$, where taking into account the assumption $\left\|x^{(k_2)}\right\|=\left\|\phi^{(k_2)}\right\|$ and the fact that $\Delta^{(k_2)}(-a_1)\geq 0$, there exist points $-a_1<t^1_{k_2+1}<t^2_{k_2+1}<t^3_{k_2+1}<0$ such that $\Delta^{(k_2+1)}(t^1_{k_2+1})<0$, $\Delta^{(k_2+1)}(t^2_{k_2+1})>0$, and $\Delta^{(k_2+1)}(t^3_{k_2+1})<0$. Continuing in the same manner, we obtain that there exist points $-a_1<t_{k_{2n-1}+1}^1<\ldots<t_{k_{2n-1}+1}^{2n}<0$ such that $(-1)^i\Delta^{(k_{2n-1}+1)}(t_{k_{2n-1}+1}^{i})>0$, $i=1,\dots,2n$, and so on all the way to the level $r-1$. 

At the level of the $(r-1)$-st derivative there exist points $-a_1<t_{r-1}^1<\ldots<t_{r-1}^{2n}<0$ such that $(-1)^i\Delta^{(r-1)}(t_{r-1}^{i})>0$, $i=1,2,\dots,2n$. In addition, $\Delta^{(r-1)}(-a_1)\geq 0$. 

This implies that on the interval $(-a_1, t^1_{r-1})$ there exists a set $S_0$ of positive measure such that $\Delta^{(r)}(t)<0$ for all $t\in S_0$. Besides that, for all $i=1,\dots, 2n-1$ on the interval $(t^i_{r-1}, t^{i+1}_{r-1})$ there exist sets $S_i\subseteq (t^i_{r-1}, t^{i+1}_{r-1})$ of positive measure such that $(-1)^i\Delta^{(r)}(t)<0$ for all $t \in S_i$. 

Therefore, the function $\Delta^{(r)}(t)$ has not fewer than $2n-1$ essential sign changes on $(-a_1, 0)$. However, this is impossible since spline $\phi$ has no more than $n$ knots and function $\Delta^{(r)}(t)$ can change sign at knots of the spline $\phi$ and no more than once on the intervals between the knots of the spline $\phi$ (and, by Lemma's assumptions, we have  $\left\|\phi^{(r)}\right\|= \left\|x^{(r)}\right\|$). The obtained contradiction proves the Lemma in the case $k_{2n}=r$.

Let now $k_{2n}<r$. Using arguments similar to the ones used above from equalities~\eqref{equalNorms} we obtain~\eqref{l2} with equality possible only if $x^{(k_1)}\equiv \phi^{(k_1)}$. Moreover, using  similar arguments once more we obtain the inequality~\eqref{l1} with equality possible only if $x^{(s)}\equiv \phi^{(s)}$. $\square$

Taking into account Theorems~\ref{AMandMomentProblem}, \ref{connection}, and Remark ~\ref{corresp}, from Theorems~\ref{momentSpaceBoundary}, \ref{momentSpaceMainRepresentation}, and \ref{momentSpaceCanonicalRepresentation} we obtain the following lemmas.

\begin{lemma}\label{boundary}
	Let integers $0 = k_1<k_2<\ldots<k_d\leq r$ be given and set $\K = (k_1,\dots,k_d)$.  Then ${M_{\K}}\in\partial A_{\K}(X)\bigcap A_{\K}(X)$  if and only if there exists $X$-perfect spline $\phi(t)$ with no more than $\frac {d-1} 2$ knots, such that
	\begin{equation}\label{*1}
		M_{\K}(\phi) = M_{\K}.
	\end{equation}
\end{lemma}

\begin{lemma}\label{mainRepresentation}
	Let integers $0 = k_1<k_2<\ldots<k_d\leq r$ be given and set $\K = (k_1,\dots,k_d)$. Assume also that ${M_{\K}}\in {\rm int}\,A_{\K}(X)$. Then there exists
	$X$-perfect spline $\phi(t)$ with  $\frac d 2$ knots such that~\eqref{*1} holds.
\end{lemma}

\begin{lemma}\label{canonicalRepresentation}
	Let positive integers $0 = k_1<k_2<\ldots<k_d\leq r$ be given and set $\K = (k_1,\dots,k_d)$. Assume also that ${M_{\K}}\in {\rm int}\,A_{\K}(X)$. Then for all $a^*>0$ there exists $X$-perfect spline $\phi(t)$ with  $\frac {d+1} 2$ knots, one of which is located at the point $a^*$, such that \eqref{*1} holds.
\end{lemma}


In addition, the following lemmas hold.

\begin{lemma}\label{l::0.0}
	Let positive integers $0<k_1<k_2<\ldots<k_{d}\leq r$ be given and set $\K = (k_1,\dots,k_{d})$. Assume that
	$M_{\K}\in \partial A_{\K}(X)\bigcap A_{\K}(X)$. 
	Then there exists $X$-perfect spline $\phi(t)$ with $s\leq\left[\frac{d - 1}{2}\right]$ knots, such that \eqref{*1} holds.
\end{lemma}

{\bf Proof.} Let function $x\in X$ be such that $M_{\K}(x) = M_{\K}$. 

First let $d=2n+1$. According to Lemma~\ref{boundary} there exists a $X$-perfect spline $\phi(t)$ with $s\leq n = \left[\frac{d-1}{2}\right]$ and a constant $C \geq 0$ such that$M_{\bf K}(\phi+C) = M_{\bf K}(x)$, where ${\bf K} = (0,k_1,\dots, k_d)$. This implies that $M_{\K}(\phi) = M_{\K}(x)$ and spline $\phi(t)$ is the required one.

Let now $d = 2n$. According to Lemma~\ref{boundary} there exists a $X$-perfect spline $\phi(t)$ with $s\leq n$ knots, such that $M_{\bf K}(\phi) = M_{\bf K}(x)$. If $s = n$, for all $C > 0$, by Lemma~\ref{boundary}, $M_{\bf K}(\phi+C)\in {\rm int} A_{\bf K}(X)$, which contradicts $M_{\K}\in \partial A_{\K}(X)\bigcap A_{\K}(X)$. Therefore, $s\leq n-1 = \left[\frac{d-1}{2}\right]$. $\square$

\begin{lemma}\label{l::0}
	Let integers $0\leq k_1<k_2<\ldots<k_{2n}\leq r$ be given and set $\K = (k_1,\dots,k_{2n})$. Let also $M_{\K}\in {\rm int} A_{\K}(X)$. 
	Then there exists $X$-perfect spline $\phi(t)$ with $n$ knots, such that~\eqref{*1} holds.
\end{lemma}
{\bf Proof.} In the case when $k_1 = 0$, the statement immediately follows from Lemma~\ref{mainRepresentation}.

Let now $k_1>0$.
Set  ${\bf K} = (0,k_1,k_2,\dots,k_{2n})$.
Let for some function $x(t)\in X$ we have $M_{\K}(x)= M_{\K}\in {\rm int} A_{\K}(X)$. Then for all $\varepsilon > 0$ $$\left(\left\|x\right\|+\varepsilon,\left\|x^{(k_1)}\right\|,\dots,\left\|x^{(k_{2n})}\right\|\right)\in {\rm int} A_{\bf K}(X).$$ Applying Lemma~\ref{mainRepresentation}, we obtain that there exists a constant $C\geq 0$ and an $X$-perfect spline $\phi(t)$ with $n$ knots, such that  $M_{\bf K}(\phi + C) = M_{\bf K}(x + \varepsilon)$. This implies that~\eqref{*1} holds. $\square$

\begin{lemma}\label{l::2}
	Let integers $0\leq k_1<k_2<\ldots<k_{2n+1}\leq r$ be given and set $\K = (k_1,\dots,k_{2n+1})$. Let also $M_{\K}\in {\rm int} A_{\K}(X)$. 
	Then for all $a^*>0$ there exists an $X$-perfect spline $\phi(t)$ with  $n+1$ knots, one of which is at $a^*$, and such that~\eqref{*1} holds.
\end{lemma}
 The proof of this Lemma is identical to the proof of Lemma~\ref{l::0}, only instead of Lemma~\ref{mainRepresentation} one should use Lemma~\ref{canonicalRepresentation}.

\begin{remark}\label{z::2}
	Under assumptions of Lemma~\ref{l::0} $($or Lemma~\ref{l::0.0}$)$, by $\phi(X,M_{\K};t)$ we denote 
	$X$-perfect spline with $d$ $\big($╨no more than $\left[\frac{d - 1}{2}\right]$, respectively$\big)$ knots, such that~\eqref{*1} holds.
\end{remark}

\begin{remark}\label{z::3}
	From Lemma~\ref{l::1} it follows that the spline $\phi(X,M_{\K};t)$ (defined in Remark~\ref{z::2}) is unique.
\end{remark}

Finally, in order to prove the main theorem we also need the following simple observation.

\begin{lemma}\label{l::-1}
	Let $\alpha>\beta>0$, $\varepsilon>0$, and function $\lambda(t)\colon\RR_+\to\RR_+$ be given. If for all sufficiently large $t>0$ $\lambda(t)\cdot t^\beta > \varepsilon$, 
	then $\lim\limits_{t\to+\infty}\lambda(t)\cdot t^\alpha=+\infty$.
\end{lemma}

\section{Solution to Kolmogorov's problem for classes $CM(\RR_-)$ and $L_{\infty,\infty}^{r,\smallsmile}(\RR_-)$}\label{S8}

Combining Theorem~\ref{conicalcone}, Lemmas~\ref{boundary} and~\ref{mainRepresentation}, as well as the connection between classes $CM(\RR_-)$ and $L_{\infty,\infty}^{r,\smallsmile}(\RR_-)$ (Theorem~\ref{connection}), we immediately obtain the solution to Kolmogorov's problem for the classes $CM(\RR_-)$ and $L_{\infty,\infty}^{r,\smallsmile}(\RR_-)$ in alternative form.

In the next theorem by $X$ we denote either class $CM(\RR_-)$ or class $L_{\infty,\infty}^{r,\smallsmile}(\RR_-)$.

\begin{theorem}\label{alternativeSolution}
	Let $d\in\NN$ and integers $0= k_1<k_2<\ldots<k_{d}\leq r$, $\K = (k_1,\dots,k_{d})$, be given. Then the set of $X$-perfect splines generates the family of admissible for $X$ sets $M_{{\bf k}}$, i.e.
	$$F_{\K}(X) = \bigg\{C+\phi(X,{\bf a},{\bf \lambda}), \; {\bf \lambda} = (\lambda_1,\dots, \lambda_{m}), \; \lambda_1,\dots, \lambda_{m} > 0,$$
	$${\bf a} = (a_1,\dots, a_{m}), \; a_1>a_2>\ldots>a_{m}>0, \;\; C\geq 0, \;\; m = \left[\frac d 2\right]\bigg\}.$$
	In the case when $d$ is even, parameter $C$ can be taken equal to zero.
\end{theorem}

Given vector ${\bf k} = (k_1,\dots,k_d)$ with integer components $0\leq k_1<k_2<\ldots<k_d\leq r$, we also introduce the following notation: ${\bf k}^2 := (k_2,\dots,k_{d})$ and $^2{\bf k}^2 := (k_2,\dots,k_{d-1})$.

\begin{theorem}\label{main}
	Let $d\in\NN$, $d\geq 3$, and nonnegative integers $0\leq k_1<k_2<\ldots<k_{d}=r$, ${\bf k} = (k_1,\dots,k_d)$, be given.  
	
	If  $d$ is odd,
	$$\{ M_{{\bf k}} \in A_{{\bf k}}(X)\}\Longleftrightarrow
		\left\{\begin{array}{c}
			M_{{\bf k}^2} \in {\rm int} A_{\K^2}(X)\\
			M_{k_1}\ge \|\phi^{(k_1)}(X,M_{{\bf k}^2})\|\\
		\end{array}\right\}
		\bigvee 
	$$
	$$
		\bigvee 
		\left\{\begin{array}{c}
			M_{{\bf k}^2} \in \partial A_{\K^2}(X)\bigcap A_{\K^2}(X)\\
			k_1>0\\
			M_{k_1}=\| \phi^{(k_1)}(X,M_{{\bf k}^2})\| \\
		\end{array}\right\}
		\bigvee
		\left\{\begin{array}{c}
			M_{{\bf k}^2} \in \partial A_{\K^2}(X)\bigcap A_{\K^2}(X)\\
			k_1=0\\
			M_{k_1}\ge \|\phi^{(k_1)}(X,M_{{\bf k}^2})\|\\
		\end{array}\right\},
	$$
	and if $d$ is even,
	$$\{ M_{{\bf k}} \in A_{{\bf k}}(X)\}\Longleftrightarrow
		\left\{\begin{array}{c}
			M_{{\bf k}^2} \in {\rm int} A_{\K^2}(X)\\
			M_{k_1} > \|\phi^{(k_1)}(X,M_{^2{\bf k}^2})\|\\
		\end{array}\right\}
		\bigvee 
	$$
	$$
		\bigvee 
		\left\{\begin{array}{c}
			M_{{\bf k}^2} \in \partial A_{\K^2}(X)\bigcap A_{\K^2}(X)\\
			k_1>0\\
			M_{k_1}=\| \phi^{(k_1)}(X,M_{^2{\bf k}^2})\| \\
		\end{array}\right\}
		\bigvee
		\left\{\begin{array}{c}
			M_{{\bf k}^2} \in \partial A_{\K^2}(X)\bigcap A_{\K^2}(X)\\
			k_1=0\\
			M_{k_1}\ge \|\phi^{(k_1)}(X,M_{^2{\bf k}^2})\|\\
		\end{array}\right\}.
	$$
	In addition, $M_{{\bf k}} \in {\rm int}A_{{\bf k}}(X)$ if and only if $M_{{\bf k}^2} \in {\rm int} A_{\K^2}(X)$ and 
	$M_{k_1}> \|\phi^{(k_1)}(X,M_{{\bf k}^2})\|$ (when $d$ is odd) or $M_{k_1}> \|\phi^{(k_1)}(X,M_{^2{\bf k}^2})\|$ (when $d$ is even).
\end{theorem}
{\bf Proof.} Due to Theorem~\ref{connection}, it is sufficient to prove the theorem in the case when $X=L_{\infty,\infty}^{r,\smallsmile}(\RR_-)$. For shortness, we denote this class simply by $X$. 

Let us first prove the necessity of the conditions. Assume that
\begin{equation}\label{MkisAdmissible}
	M_{{\bf k}} \in A_{{\bf k}}. 
\end{equation}
Let a function $x\in L^{r,\smallsmile}_{\infty,\infty}({\RR_-})$ be such that $M_{\bf k}(x)=M_{\bf k}$.

Since the set $A_{\K^2}$ is convex, then $A_{\K^2} = {\rm int}A_{\K^2} \bigcup \left(\partial A_{\K^2}\bigcap A_{\K^2}\right)$. It implies that $M_{\K^2}\in {\rm int}A_{\K^2}$ or $M_{\K^2}\in \partial A_{\K^2}\bigcap A_{\K^2}$.  Let first
\begin{equation}\label{Mk2IsInternal}
	M_{{\bf k}^2} \in {\rm int} A_{\K^2}.
\end{equation}
The necessity of the condition 
\begin{equation}\label{***}
	M_{k_1}\ge \|\phi^{(k_1)}(M_{{\bf k}^2})\|
\end{equation}
in the case when $d=2n+1$ follows from Lemma~\ref{l::1}. Indeed, 
$$\|\phi^{(k_i)}(M_{{\bf k}^2})\|=\|x^{(k_i)}\|, \qquad i=2,3,\dots, 2n+1.$$
 Moreover, from Lemma~\ref{l::0} we obtain that the spline $\phi(M_{{\bf k}^2})$ has $n$ knots. From Lemma~\ref{l::1} it now follows that $$M_{k_1}=\|x^{(k_1)}\|\geq \|\phi^{(k_1)}(M_{{\bf k}^2})\|.$$
Let now $d=2n$. Then $$\|\phi^{(k_i)}(M_{^2{\bf k}^2})\|=\|x^{(k_i)}\|, \qquad i=2,3,\dots, 2n-1,$$ and the spline $\phi(M_{^2{\bf k}^2})$ has $n-1$ knots due to Lemma~\ref{l::0}. From Lemma~\ref{l::1} we have $M_{k_1}=\|x^{(k_1)}\|\geq \|\phi^{(k_1)}(M_{^2{\bf k}^2})\|$ with equality possible only if $x^{(k_1)}\equiv \phi^{(k_1)}(M_{^2{\bf k}^2})$.  Lemma~\ref{l::1} implies that for all $\varepsilon>0$  $(M_{k_2},\dots, M_{k_{2n-1}}, \|\phi^{(r)}(M_{^2{\bf k}^2})\|-\varepsilon)\notin A_{{\bf k}^2}$.  Hence, from~\eqref{Mk2IsInternal} we conclude that $\|\phi^{(k_d)}(M_{^2{\bf k}^2})\|<M_{k_d}$, which, in particular, implies 
\begin{equation}\label{even}
	M_{k_1} > \|\phi^{(k_1)}(M_{^2{\bf k}^2})\|.
\end{equation}
Next we consider the case when
\begin{equation}\label{Mk2isOnBoundary}
	M_{{\bf k}^2} \in \partial A_{\K^2}(X)\bigcap A_{\K^2}(X).
\end{equation}
Let us first take $d = 2n+1$. In this case the spline $\phi(M_{\K^2})$ has no more than $n-1$ knots due to Lemma~\ref{l::0.0}. Moreover, $$\|\phi^{(k_i)}(M_{{\bf k}^2})\|=\|x^{(k_i)}\|, \qquad i=4,\dots, 2n+1$$ and $\|\phi^{(k_2)}(M_{{\bf k}^2})\|=\|x^{(k_2)}\|$. Then Lemma~\ref{l::1}  implies that $\phi^{(k_2)}(M_{{\bf k}^2})\equiv x^{(k_2)}$. If $k_1>0$, then $\phi^{(k_1)}(M_{{\bf k}^2})\equiv x^{(k_1)}$ and, hence, $\|\phi^{(k_1)}(M_{{\bf k}^2})\|= M_{k_1}.$ If $k_1=0$ then $x\equiv \phi(M_{{\bf k}^2}) + C,$ $C\geq 0$, and, hence, $\|\phi(M_{{\bf k}^2})\|\leq M_{k_1}.$

The case of even $d$ can be considered in a similar way.
 


In the case when $d$ is odd, ~\eqref{Mk2IsInternal} holds, and \eqref{***} holds with equality, we have
\begin{equation}\label{MkisOnBoundary}
	M_{\K}\in \partial A_{\K}\bigcap A_{\K},
\end{equation}
since from Lemma~\ref{l::1} it follows that for all $\varepsilon > 0$ it is true that $(M_{k_1}-\varepsilon,M_{k_2},\dots,M_{k_{d}})\notin A_{\K}$. 

This concludes the proof of necessity of the assumptions of the Theorem.

Next we prove that they are also sufficient.

First, we observe that in the case when~\eqref{Mk2isOnBoundary} holds and $d$ is odd, the sufficiency of the Theorem's assumptions is obvious. In the case of even $d=2n$ we note that from~\eqref{Mk2isOnBoundary} and Lemma~\ref{l::1} it follows that $\|\phi^{(r)}(M_{^2{\bf k}^2})\| = M_{k_d}$. Indeed, let function $x\in L^{r,\smallsmile}_{\infty,\infty}({\RR_-})$ be such that $M_{{\bf k}^2}(x)=M_{{\bf k}^2}$. The spline $\phi(M_{^2{\K^2}})$ has no more than  $n-2$ knots due to Lemma~\ref{l::0.0}. Moreover,  $\|\phi^{(k_i)}(M_{^2{\bf k}^2})\|=\|x^{(k_i)}\|$, $i=4,5,\dots, 2n-1$, and $\|\phi^{(k_2)}(M_{^2{\bf k}^2})\|=\|x^{(k_2)}\|$. Lemma~\ref{l::1} now implies that $\phi^{(k_2)}(M_{^2{\bf k}^2})\equiv x^{(k_2)}$, and, hence, $\phi^{(r)}(M_{^2{\bf k}^2})\equiv x^{(r)}$. This means that $\|\phi^{(r)}(M_{^2{\bf k}^2})\|\equiv \|x^{(r)}\| = M_{k_d}$.

Next we prove the sufficiency in the case when~\eqref{Mk2IsInternal} is satisfied. 

Let us consider the case of odd $d=2n+1$. The sufficiency is obvious in the case $k_1 = 0$, and, hence, we may assume $k_1 > 0$.

Without loss of generality we may assume that~\eqref{***} holds with strict inequality.  
From~\eqref{Mk2IsInternal} it follows that for any $\varepsilon > 0$ $(\|\phi(M_{\K^2})\|+\varepsilon,M_{k_2},\dots, M_{k_d})\in {\rm int}A_{{\bf K}}$, where ${\bf K} = (0,k_2,k_3,\dots, k_d)$. Let $a_n$ be the rightmost knot of the spline $\phi(M_{\K^2})$. By Lemma~\ref{l::2} for all $\varepsilon > 0$ there exists $X$-perfect spline $\psi = \psi(\varepsilon)=\phi({\bf a}(\varepsilon),{\bf \lambda}(\varepsilon))$, where the vectors ${\bf a}(\varepsilon)=(a_1(\varepsilon),a_2(\varepsilon),\dots,a_{n+1}(\varepsilon))$ and ${\bf \lambda(\varepsilon) }= (\lambda_1(\varepsilon),\dots,\lambda_{n+1}(\varepsilon))$ are such that  $a_1(\varepsilon)>a_2(\varepsilon)>\ldots>a_{n+1}(\varepsilon)= \frac {a_n}{2}$, $\lambda_1(\varepsilon),\dots,\lambda_{n+1}(\varepsilon)>0$, for $i=2,\dots, d$ $\left\|\psi^{(k_i)}(\varepsilon)\right\| = M_{k_i}$ and $\left\|\psi(\varepsilon)\right\| = \varepsilon + \left\|\phi(M_{\K^2})\right\|$. 

For $i=1,\dots,n+1$ we set $a_i^*:=\lim\limits_{\varepsilon\to+\infty}a_i(\varepsilon)$ (limit can be finite or infinite). It is clear that $a_1^*=\infty$. Let us assume that the function $\lambda_1(\varepsilon)a_1^{r-k_1}(\varepsilon)$ is bounded as  $\varepsilon \to+\infty$. Then by Lemma~\ref{l::-1}, we have $\lambda_1(\varepsilon)a_1^{r-k_i}(\varepsilon)\to 0$ as $\varepsilon \to+\infty$ for all $i=2,\dots,d$. This implies that $a_2^*=\infty$, since by Remark~\ref{z::3} and due to the fact that $a_{n+1}(\varepsilon)\equiv \frac {a_n}{2}$, all limits $a_i^*$, $i=2,\dots, n+1$, cannot be finite. Repeating similar arguments, we obtain that for some $i=1,\dots,n$ the function $\lambda_i(\varepsilon)a_i^{r-k_1}(\varepsilon)$ is unbounded as $\varepsilon\to\infty$. It means that the value of the norm $\left\|\psi^{(k_1)}(\varepsilon)\right\|$ can be made arbitrarily large. Due to convexity of the set $A_{{\bf k}}$, it implies that~\eqref{MkisAdmissible} holds.

In addition, since the fact that~\eqref{Mk2IsInternal} holds for an arbitrary $M_{k_1}$, for which strict inequality~\eqref{***} holds, implies~\eqref{MkisAdmissible}, then for such $M_{k_1}$ we have $M_{\K}\in {\rm int}A_{{\bf k}}$. 

Finally, let $d = 2n$. 
By Lemma~\ref{l::2}, for all $a>0$ there exists spline $\psi(a)$ with $n$ knots, one of which is at point $a$, such that  we have $M_{\K^2}(\psi(a)) = M_{\K^2}$. As $a\to 0$ we have pointwise convergence $\psi(a)\to\phi(M_{^2{\bf k}^2})$, and, hence, for small enough $a$ we have $\left\|\psi^{(k_1)}(a)\right\|<M_{k_1}$. Taking the limit as $a\to\infty$ and using arguments similar to the case of odd $d$, we obtain that~\eqref{MkisAdmissible} holds. $\square$

\section{Some applications}\label{S9}
\subsection{On the smoothest Hermite-Birkhoff interpolation}

We now consider the following problem. One has to find the necessary and sufficient conditions to guarantee the existence of a function $f\in X$ such that
\begin{equation}\label{interpolation}
	f^{(k_i)}(0) = M_{k_i},\qquad i =1,\dots,2n,
\end{equation}
where $0\leq k_1<k_2<\ldots<k_{2n}<r$, $n\in\NN$. 
This problem is called Hermite-Birkhoff interpolation at $0$.  If the problem above has a solution, then it is a problem of interest to find a function, which satisfies conditions~\eqref{interpolation} and has the least possible value of the norm of the highest derivative. We call this problem a problem of smoothest Hermite-Birkhoff interpolation. For more information about the smoothest interpolation see, for instance, works ~\cite{Holladay, Karlin75, Pinkus}, and references therein. For investigations on Hermite-Birkhoff interpolation we refer the reader to works \cite{Lorentz, LorentzZeller, KarlinKaron, Schoenberg}, and references therein.

It is clear that the problem of existence of interpolation function is equivalent to Kolmogorov's problem (as uniform norm of  monotone function is attained at point $0$). Moreover, from the results of the present paper (see Lemma~\ref{l::1}) it follows that the spline $\phi(X, M_{\bf k})$ has the lowest value of the  {uniform} norm of highest derivative among all functions from $X$ satisfying~\eqref{interpolation}, $M_{\bf k} = (M_{k_1}, \dots, M_{k_d})$.

\subsection{On sharp estimates for intermediate moments}
Let $0\leq k_1<\ldots<k_{2m}\leq r$ and $u_s(t)=t^{k_s}$, $t\geq 0$, $s=1,\dots,2m$. We consider the following problem. On the set of solutions to moment problem~\eqref{moments} (with $n=2m$), one needs to find the best estimates for the moments
$$
c_p= \int_{0}^\infty t^p d\sigma(t), \;\;\; p\neq k_s,\;\;\; s=1,\dots,2m.
$$
This problem belongs to the class of problems on finding extremal values of integrals, which were at first studied in the works of Chebyshev and Markov (see, for instance, ~\cite{Cheb, Kr, Markov}).

Using similar to the proof of Lemma~\ref{l::1} arguments one can prove the following theorem.

\begin{theorem}\label{l1Generalization1}
Let $x\in X$, $k_0:=0 \leq k_1<\ldots<k_{2m}\leq r$. Let ${X}$-spline $\phi(t)$ with no more than $m$ knots be such that 
\begin{equation*}
\left\|\phi^{(k_i)}\right\| = \left\|x^{(k_i)}\right\|,\qquad i = 1,\dots,2m.
\end{equation*}
If for some $s\in\{1,\dots, 2m-1\}$ $k_s<p<k_{s+1}$, $k_{s}\leq p<k_{1}$ with $s=0$, or  $k_{s}<p\leq r$ with $s=2m$,  then
$$
		(-1)^s\left\|\phi^{(p)}\right\|\leq (-1)^s\left\|x^{(p)}\right\|.
$$
\end{theorem}

Theorem~\ref{l1Generalization1}, together with the connection between moment problem and Kolmogorov's problem on the class of absolute monotone functions, gives sharp estimates (from above or from below depending on $p$) for the intermediate moment $c_p$.

\subsection{On extremal distribution functions}
We define the set $\mathcal{F}$ of functions $F\in X$ as follows
$$
F^{(k_i)}(0)=M_{k_i},\qquad  i=1,\dots,2n,
$$
$$F(-\infty)=0,$$
where $0 = k_1<\ldots<k_{2n}\leq r$, and $M_{k_1}=1$. Such functions can be considered as distribution functions for random variables that take nonpositive values. 

Let $A<0$ be given. The problem is to find the distribution function $F\in\mathcal{F}$ such that the corresponding random variable $\xi=\xi_F$ has the highest probability $P(\xi >A)$.

Using similar to the proof of Lemma~\ref{l::1} arguments one can prove the following theorem.

\begin{theorem}\label{l1Generalization2}
Let $x\in X$, $0 \leq k_1<\ldots<k_{2n}\leq r$. Let ${X}$-spline $\phi(t)$ with no more than $n$ knots be such that 
\begin{equation*}
\left\|\phi^{(k_i)}\right\| = \left\|x^{(k_i)}\right\|,\qquad i = 1,\dots,2n.
\end{equation*}
Then for all $t\leq 0$
$$\phi^{(k_1)}(t)\leq x^{(k_1)}(t).$$
\end{theorem}

From Theorem~\ref{alternativeSolution} it follows that either $\mathcal{F}$ is empty, or there exists a spline $\phi\in \mathcal{F}$ with no more than $n$ knots. Then $\phi$ is the extremal distribution function due to Theorem~\ref{l1Generalization2}.

\section{Acknowledgements}
This project was supported by Simons Collaboration Grant N. 210363.

\end{document}